\documentclass[12pt,oneside,english]{amsart}

\usepackage[english]{babel}
\usepackage{multirow}
\usepackage{graphicx}
\usepackage{amssymb}
\usepackage{yhmath}
\usepackage{verbatim} 
\usepackage{amsmath}
\usepackage{amsthm}
\usepackage{color}
\usepackage{epstopdf} 
\usepackage{wrapfig}
\usepackage[all]{xy}
\usepackage{hyperref}
\usepackage{enumitem} 

\newcommand{\BE}{\begin{equation}}
\newcommand{\EE}{\end{equation}}

\newcommand{\beq}{\begin{eqnarray}}
\newcommand{\eeq}{\end{eqnarray}}
\newcommand{\nn}{\nonumber}
\newcommand{\bev}{\begin{verbatim}}
\newcommand{\eev}{\end{verbatim}}

\newtheorem{theorem}{Theorem}
\newtheorem{coro}{Corollary}

\newcounter{example}[section]

\calclayout
\makeatletter

\begin{document}
\title{Degree of the Gauss map and curvature integrals for closed hypersurfaces}
\author{Fabiano G. B. Brito}

\author{Icaro Gon\c calves}

\address{Centro de Matem\'atica, Computa\c{c}\~ao e Cogni\c{c}\~ao,
Universidade Federal do ABC,
09.210-170 Santo Andr\'e, Brazil}
\email{fabiano.brito@ufabc.edu.br}

\address{Dpto. de Matem\'{a}tica, Instituto de Matem\'{a}tica e Estat\'{i}stica,
Universidade de S\={a}o Paulo, R. do Mat\={a}o 1010, S\={a}o Paulo-SP
05508-900, Brazil.}
\email{icarog@ime.usp.br}

\subjclass[2010]{57R25, 57R30, 47H11}

\thanks{The second author was supported by CNPq, 141113/2013-8, Brasil}

\begin{abstract}
Given a unit vector field on a closed Euclidean hypersurface, we define a map from the hypersurface to a sphere in the Euclidean space. This application allows us to exhibit a list of topological invariants which combines the second fundamental form of the hypersurface and the vector field itself. We show how these invariants can be used as obstructions to certain codimension one foliations. 
\end{abstract}
\maketitle

\section{Introduction}
Let $M$ be a closed Riemannian manifold of constant sectional curvature, and let $\mathcal{F}$ be a transversely oriented codimension one foliation on $M$. 
Brito et al, \cite{brito}, have shown that the total higher-order mean curvatures of the leaves, i.e. the integral of the elementary symmetric functions of the eigenvalues of the second fundamental form of the leaves, do not depend on the foliation. 
In that direction, these integral formulas for codimension one foliations have been generalized by means of Newton transformations and for any closed Riemannian manifold, see \cite{pawel1} and \cite{pawel2} and the references therein. 

The total higher-order mean curvatures were computed in \cite{brito} through an application defined by Milnor in \cite{milnor1}; for the case of positive curvature, let $\varphi_t: \mathbb{S}^{n+1}\to \mathbb{S}^{n+1}(\sqrt{1+t^2})$, $\varphi_t(x) = x + tN(x)$, where $N$ is the vector field normal to $\mathcal{F}$. The integrals come from the determinant of the Jacobian matrix of $\varphi_t$, and coefficient comparison in a polynomial in the variable $t$. 

All the known cases, for example \cite{brito}, \cite{pawel1} and \cite{pawel2}, the integral formulae show independence from the foliation, but they still depend on the geometry of the ambient manifold. In this note we define a map similar to $\varphi_t$, in which we employ a unit vector field on a closed Euclidean hypersurface $M^{n+1}$ immersed in $\mathbb{R}^{n+2}$, and also the normal map of the hypersurface. This map provides integral formulae, in which we can recover the aforementioned positive case, but they rather depend on the topology of $M$, and this dependence comes through the degree of the normal map of $M$. 

\section{The map $\varphi^{\vec{v}}_{t}$}
Let $M^{n+1}$ be an oriented closed immersed Euclidean hypersurface, with $\chi(M) = 0$, and let $\vec{v}: M \to TM$ be a smooth unit vector field on $M$. We fix an orientation of $M$ so the normal map $\nu : M^{n+1}\to \mathbb{S}^{n+1}$ is well defined, $\nu(x) = N(x)$, where $N$ is a unitary field, normal to $M$. 

For a real number $t>0$, define $\varphi^{\vec{v}}_{t}:  M^{n+1}\to \mathbb{R}^{n+2}$, by $\varphi^{\vec{v}}_{t}(x) = \nu(x) + t\vec{v}(x)$. Both $\nu(x)$ and $\vec{v}(x)$ are unitary in $\mathbb{R}^{n+2}$, for all $x \in M$; thus, $||\varphi^{\vec{v}}_{t}(x)||^2 = 1+t^2$, which means that the image of $\varphi^{\vec{v}}_{t}$ lies in the sphere $\mathbb{S}^{n+1}(\sqrt{1+t^2})$ of radius $\sqrt{1+t^2}$. On the other hand, $\nu$ and $\vec{v}$ are smooth maps, and adding is a smooth operation in $\mathbb{R}^{n+2}$, so $\varphi^{\vec{v}}_{t}$ is again a smooth application. 

The degree formula (see for example  Guillermin-Pollack, page 188) says that for any smooth map $f:X\to Y$ between manifolds with the same dimension $k$, and for a $k$-form $\omega$ on $Y$, 
\begin{equation}
\label{degree}
\int_{X}f^* \omega = \deg(f)\int_Y \omega 
\end{equation}

In the case of $\varphi^{\vec{v}}_{t} :  M^{n+1}\to \mathbb{S}^{n+1}(\sqrt{1+t^2})$, set $\omega$ as the volume form of $\mathbb{S}^{n+1}(\sqrt{1+t^2})$, and ${\rm vol}(\mathbb{S}^{n+1})$ as the volume of the unit Euclidean sphere. The RHS of (\ref{degree}) reads
\begin{eqnarray}
\label{RHS}
\deg(\varphi^{\vec{v}}_{t})\int_{\mathbb{S}^{n+1}(\sqrt{1+t^2})} \omega &=& \deg(\varphi^{\vec{v}}_{t}){\rm vol}(\mathbb{S}^{n+1}(\sqrt{1+t^2}))\nn\\ &=& \deg(\varphi^{\vec{v}}_{t}) {\rm vol}(\mathbb{S}^{n+1})(\sqrt{1+t^2})^{n+1}\nn\\ &=& \deg(\varphi^{\vec{v}}_{t}) {\rm vol}(\mathbb{S}^{n+1})(\sqrt{1+t^2})(\sqrt{1+t^2})^n
\end{eqnarray}

For the LHS of (\ref{degree}), there exists a relation between the pullback of $\omega$ via $\varphi^{\vec{v}}_{t}$ and the volume form of $M$, $\omega_M$; it follows, 
\begin{equation}
\label{LHS}
\int_{M}(\varphi^{\vec{v}}_{t})^* \omega = \int_M \det(d \varphi^{\vec{v}}_{t}) \omega _M
\end{equation}

\section{Topological invariants}
Since the determinant is invariant under similarity transformations, we may choose any orthonormal basis for $T_x M^{n+1}$ and $T_{\varphi^{\vec{v}}_{t}(x)}\mathbb{S}^{n+1}(\sqrt{1+t^2})$ in order to compute the determinant of $d(\varphi^{\vec{v}}_{t})$. For $T_x M^{n+1}$, take the following adapted set of orthonormal vectors
\begin{equation}
\lbrace e_1, \dots, e_{n}, \vec{v} \rbrace.\nn
\end{equation}
We notice that $\varphi^{\vec{v}}_{t}(x) \perp e_i$, for every $1\leq i\leq n$, and that $\varphi^{\vec{v}}_{t}(x)$ can be chosen as normal to $\mathbb{S}^{n+1}(\sqrt{1+t^2})$. So we may complete the set $\{e_1, \dots, e_{n}\}$ with a vector which is simultaneously normal to all $e_i$, $1\leq i\leq n$, and $\varphi^{\vec{v}}_{t}(x)$. Thus
\begin{equation}
\left\lbrace e_1, \dots, e_{n}, u:= \frac{\vec{v}}{\sqrt{1+t^2}} - t\frac{N}{\sqrt{1+t^2}}\right\rbrace\nn
\end{equation}
is an orthonormal basis for $T_{\varphi^{\vec{v}}_{t}(x)}\mathbb{S}^{n+1}(\sqrt{1+t^2})$.

If $S = S_N: TM \to TM$ is the shape operator of $M$, then $d\nu(\cdot) = S(\cdot) = (D_{(\cdot)} N)^{\top}$, where $(\cdot)^{\top}$ is the projection on $TM$, and $D$ denotes the Euclidean connection. Setting $\nabla$ as the induced Levi-Civita connection of $M$, we have
\begin{equation}
\label{dv}
d\vec{v}(\cdot) = D_{(\cdot)}\vec{v} = \nabla_{(\cdot)}\vec{v} + (D_{(\cdot)}\vec{v})^{\perp} = \nabla_{(\cdot)}\vec{v} + \left\langle D_{(\cdot)}\vec{v}, N\right\rangle N = \nabla_{(\cdot)}\vec{v} - \left\langle \vec{v}, S(\cdot)\right\rangle N.
\end{equation}

With respect to the aforementioned notation the acceleration $\nabla_{\vec{v}}\vec{v}$ has the following components: $v_i :=\left\langle \nabla_{\vec{v}}\vec{v}, e_i\right\rangle$. Evidently, $v_{n+1} = \left\langle \nabla_{\vec{v}}\vec{v}, \vec{v}\right\rangle = 0$. In addition, we define $a_{ij} := \left\langle \nabla_{e_i}\vec{v}, e_j\right\rangle$ and $h_{AB} := \langle S(e_A), e_B\rangle$, for $1\leq i,j \leq n$ and $1\leq A, B \leq n+1$; $h_{AB}$ are the entries of the second fundamental form of $M$. The matrix $(a_{ij})$ locally describes the behavior of the normal bundle of $\vec{v}$; e.g. $(a_{ij})$ is symmetric if and only if the normal bundle is integrable to a codimension one foliation of $M$. From the equation (\ref{dv}), $d\vec{v}(e_i) = \nabla_{e_i}\vec{v} - \langle \vec{v}, S(e_i)\rangle N = \nabla_{e_i}\vec{v} - h_{n+1\;i} N,
$ and $d\vec{v}(\vec{v}) = \nabla_{\vec{v}}\vec{v} - \langle \vec{v}, S(\vec{v})\rangle N = \nabla_{\vec{v}}\vec{v} - h_{n+1\ n+1} N$.

We have
\begin{eqnarray}
\left\langle d\varphi^{\vec{v}}_{t}(e_i), e_j\right\rangle &=& \left\langle d\nu(e_i) + td\vec{v}(e_i), e_j \right\rangle = \left\langle S(e_i) + t\nabla_{e_i}\vec{v} - th_{n+1\;i} N, e_j\right\rangle = h_{ij} + ta_{ij}\nn\\ 
\nn\\
\left\langle d\varphi^{\vec{v}}_{t}(e_i), u \right\rangle &=& \left\langle S(e_i) + t\nabla_{e_i}\vec{v} - th_{n+1\;i} N, \frac{\vec{v}}{\sqrt{1+t^2}} - t\frac{N}{\sqrt{1+t^2}}\right\rangle\nn\\ 
&=& \left\langle S(e_i), \frac{\vec{v}}{\sqrt{1+t^2}}\right\rangle - \left\langle S(e_i), t\frac{N}{\sqrt{1+t^2}}\right\rangle
+ \left\langle t\nabla_{e_i}\vec{v}, \frac{\vec{v}}{\sqrt{1+t^2}}\right\rangle\nn\\ 
&-& \left\langle t\nabla_{e_i}\vec{v}, t\frac{N}{\sqrt{1+t^2}}\right\rangle - \left\langle th_{n+1\;i} N, \frac{\vec{v}}{\sqrt{1+t^2}}\right\rangle + 
\left\langle th_{n+1\;i} N, t\frac{N}{\sqrt{1+t^2}}\right\rangle\nn\\
&=& \frac{1}{\sqrt{1+t^2}}\left\langle S(e_i), \vec{v}\right\rangle + \frac{t}{\sqrt{1+t^2}}\left\langle \nabla_{e_i}\vec{v}, \vec{v}\right\rangle + \frac{t^2}{\sqrt{1+t^2}}h_{n+1\;i} \left\langle  N, N \right\rangle\nn\\
&=& \frac{1}{\sqrt{1+t^2}}h_{i\; n+1} + \frac{t^2}{\sqrt{1+t^2}}h_{n+1\;i}\nn\\
&=& \sqrt{1+t^2}h_{n+1\;i}\nn\\
\nn\\
\left\langle d\varphi^{\vec{v}}_{t}(\vec{v}), e_i \right\rangle &=& \left\langle d\nu(\vec{v}) + td\vec{v}(\vec{v}), e_i \right\rangle = \left\langle S(\vec{v}) + t\nabla_{\vec{v}}\vec{v} - t h_{n+1\; n+1} N, e_i \right\rangle = h_{n+1\;i} + tv_i\nn\\
\nn\\
\left\langle d\varphi^{\vec{v}}_{t}(\vec{v}), u \right\rangle &=& \left\langle d\varphi^{\vec{v}}_{t}(\vec{v}), \frac{\vec{v}}{\sqrt{1+t^2}} - t\frac{N}{\sqrt{1+t^2}}\right\rangle \nn\\ 
&=& \left\langle S(\vec{v}) + t\nabla_{\vec{v}}\vec{v} - t h_{n+1\; n+1} N, \frac{\vec{v}}{\sqrt{1+t^2}} - t\frac{N}{\sqrt{1+t^2}}\right\rangle \nn\\
&=& \frac{1}{\sqrt{1+t^2}} h_{n+1\; n+1}  + \frac{t^2}{\sqrt{1+t^2}} h_{n+1\; n+1}\nn\\ 
&=& \sqrt{1+t^2}h_{n+1\; n+1}\nn
\end{eqnarray}

Then
\[
d\varphi^{\vec{v}}_{t} = 
  \left(\begin{array}{ccc|c}
    \multicolumn{3}{c|}{\multirow{3}{*}{\raisebox{2pt}{$h_{ij} + ta_{ij}$}}}     & \sqrt{1+t^2}h_{n+1\;1}           \\ 
    &       &            & {\vdots}    \\
    &       &            & \sqrt{1+t^2}h_{n+1\; n}           \\ \hline
    h_{n+1\;1} + tv_1       & \cdots & h_{n+1\;n} + tv_{n} & \sqrt{1+t^2}h_{n+1\; n+1}
  \end{array}\right)
\]

Define some column vectors by 
$$
V_j = (a_{1j},\dots, a_{n\;j}, v_j), \quad H_{j} = (h_{1j}, \dots, h_{n+1\;j}), \quad H_{n+1} = (h_{n+1\; 1}, \dots, h_{n+1\;n}, h_{n+1\; n+1}),
$$
and rewrite the above matrix as follows
\begin{equation}
\label{dvarphi-vector}
d\varphi^{\vec{v}}_{t} = 
	\left(
	\begin{array}{cccc} 
	H_1 + tV_1 & \cdots & H_{n} + tV_{n} & \sqrt{1+t^2}H_{n+1}
	\end{array}
	\right).
\end{equation}

The determinant is linear with respect to the sum of column vectors, so the equation (\ref{dvarphi-vector}) simplify computations concerning an explicit formula for $\det(d\varphi^{\vec{v}}_{t})$ written in terms of the second fundamental form of $M$ and the components depending on the normal bundle of $\vec{v}$. Therefore,
\begin{eqnarray}
\det(d\varphi^{\vec{v}}_{t}) &=& \sqrt{1+t^2} \sum_{k=0}^{n}\eta_k t^k, \quad {\rm where}\nn\\
                      \eta_0 &=& \det(h_{AB})\nn\\
                      \eta_1 &=& \sum_i \det\left(\begin{array}{ccccc}H_1 & \cdots & V_i & \cdots & H_{n+1}\end{array}\right)\nn\\
                      \eta_2 &=& \sum_{i<j} \det\left(\begin{array}{ccccccc}H_1 & \cdots & V_i & \cdots & V_j & \cdots & H_{n+1}\end{array}\right)\nn\\
                      	    &\vdots & \nn\\
                   \eta_{n} &=&	\det\left(\begin{array}{cccc}V_1 & \cdots & V_{n} & H_{n+1}\end{array}\right) \nn
\end{eqnarray}

From the equations (\ref{RHS}) and (\ref{LHS}), 
\begin{equation}
\sum_{k=0}^{n} t^k \int_{M} \eta_k =
\begin{cases}
\deg(\varphi^{\vec{v}}_{t}) {\rm vol}(\mathbb{S}^{n+1})\sum_{k=0}^{n/2}{n/2 \choose k} t^{2k},
  &\mbox{if}\;n \; \mbox{is}\;\mbox{even},\\
\deg(\varphi^{\vec{v}}_{t}) {\rm vol}(\mathbb{S}^{n+1})\sqrt{1+t^2}\sum_{k=0}^{(n-1)/2}{(n-1)/2 \choose k} t^{2k}, & \mbox{if}\; n \; \mbox{is}\;\mbox{odd}.
\end{cases}
\nn
\end{equation}

We first observe that when $n$ is even both sides of the last equation are polynomials in $t$, and the powers of the RHS are all even integers. This implies that the coefficients multiplying odd powers on the LHS are all zero. When $n$ is odd, the appearence of $\sqrt{1+t^2}$ shows that all coefficients on the LHS are zero.  On the other hand, for $t$ sufficiently small we have that the degrees of the maps $\varphi^{\vec{v}}_{t}$ and $\nu$ coincide, $\deg(\varphi^{\vec{v}}_{t}) = \deg(\nu)$. So we conclude 
\begin{equation}\boxed{
\label{integral-formula}
\int_{M} \eta_k=
\begin{cases}
\deg(\nu){n/2 \choose k/2} {\rm vol}(\mathbb{S}^{n+1}),
  &\mbox{if}\; k \;\mbox{and}\; n \; \mbox{are}\;\mbox{even},\\
0, & \mbox{if}\; k \;\mbox{or}\; n \; \mbox{is}\;\mbox{odd}.
\end{cases}}
\end{equation}

\section{The degree of the normal map}
Set $\deg (\nu) = d$. In \cite{milnor2}, we can find an answer to the question of whether $d$ assumes any given integer,
\begin{theorem}[Milnor, \cite{milnor2}]
\label{thm-milnor-1}
For $n$ odd,

(a) If M can be immersed in $\mathbb{R}^{n+1}$ with degree $d$, then it can also be immersed with any degree $d'$ which is congruent to $d$ modulo $2$.

(b) If $M$ can be immersed in $\mathbb{R}^{n+1}$ with normal degree zero, then $M$ is parallelizable.

(c) Assume $M$ is not parallelizable. If $M$ can be immersed in $\mathbb{R}^{n+1}$ at all, then it can be immersed with arbitrary odd degree, but cannot be immersed with even degree.
\end{theorem}

The first item asserts how plentiful is the set of possible values for $d$. Nevertheless, when the Betti numbers of $M$ are taken into account, restrictions arise by comparison to the normal degree,

\begin{theorem}[Milnor, \cite{milnor2}]
\label{thm-milnor-2}
For any imbedding of $M^n$ in $\mathbb{R}^{n+1}$ the degree $d$ of the normal map satisfies
\begin{equation}
2d \equiv \beta(M)\mod{2}, \quad 2|d|\leq \beta(M),
\end{equation}
where $\beta(M) = \beta_0(M) + \beta_1(M) + \cdots + \beta_n(M)$ denote the sum of the Betti numbers of the manifold $M$. In addition, if $M$ is oriented, then $2-\frac{1}{2}\beta(M)\leq d \leq \frac{1}{2}\beta(M)$.
\end{theorem}

For a number $d$ occur as the degree of a immersion, the ring $H^{*}(M)$ must split into the sum of two subalgebras, and
$d$ should appear as the Euler characteristic of one of the two. For further discussion, e.g. proofs of the theorems as well as examples, see \cite{milnor2}. See also \cite{smale} for the case of possible values of the degree for immersion of spheres.

\section{First consequences of (\ref{integral-formula})}
In this section we summarize some direct consequences of the relation \ref{integral-formula}. 
\subsection{Codimension one foliations}

Let $\mathcal{F}^{2n}$ be a transversely oriented codimension one foliation of a closed oriented immersed hypersurface $M^{2n+1}$ of $\mathbb{R}^{2n+2}$. Take $\vec{v}$ as a unit vector field normal to the leaves. In this case, at each point of $M$, the matrix $(a_{ij})$ represent the second fundamental form of the leaf passing through that point, and thus it is symmetric. 
Then, 
$$
\eta_{2n} = 
  \det \left(\begin{array}{ccc|c}
    \multicolumn{3}{c|}{\multirow{3}{*}{\raisebox{2pt}{$a_{ij}$}}}     & h_{2n+1\;1}           \\ 
    &       &            & {\vdots}    \\
    &       &            & h_{2n+1\; 2n}           \\ \hline
    v_1       & \cdots & v_{2n} & h_{2n+1\;2n+1}
  \end{array}\right). 
$$
If $\eta_{2n} = 0$ (or for any other $1\leq k \leq n$, $\eta_{2k} = 0$) then (\ref{integral-formula}) implies that $\deg(\nu) = 0$. 
\begin{theorem}
If $M^{2n+1}$ admits a transversely oriented codimension one foliation, such that the second fundamental form of the leaves has rank less or equal than $2(n-1)$, then $\deg(\nu) = 0$, where $\nu: M^{2n+1}\to \mathbb{S}^{2n+1}$ is the normal map. 
\end{theorem}
Combining last theorem together with item (b) of \ref{thm-milnor-1} imply that $M^{2n+1}$ is parallelizable. On the other hand, a totally geodesic foliation of $M$ clearly satisfies the hypothesis on the rank of the second fundamental form. Thus,  

\begin{coro}
If $M$ has a codimension one totally geodesic foliation, then the degree of its normal map must be zero.  
\end{coro}

Taking theorem \ref{thm-milnor-1} into account, we may conclude that between all immersions of a given closed hypersurface, and all possible values $\deg(\nu)$ can assume, foliations satisfying ${\rm rank}(a_{ij})\leq 2(n-1)$ (and therefore totally geodesic ones) might occur only when $\deg(\nu) = 0$. In the spirit of itens $(b)$ and $(c)$ of theorem \ref{thm-milnor-1}, Smale, in \cite{smale}, proved that there exists an immersion of a sphere in the Euclidean space with normal degree zero, if such a sphere is parallelizable. Thus, the only possible spheres that could admit a codimension one totally geodesic foliation are $\mathbb{S}^3$ and $\mathbb{S}^7$ immersed with degree zero. By Sullivan's theorem, \cite{sullivan}, $\mathbb{S}^3$ does not admit such a foliation for any Riemannian metric; one just has to use the existence of a Reeb component. Therefore, the following question remains: is there an immersion of $\mathbb{S}^7$ in the Euclidean space endowed with a codimension one totally geodesic foliation? It would be interesting to know whether Ghys' characterization of those foliations in \cite{ghys} provides a quick answer to this question.

\section{Conclusion and further research}
In a forthcoming version of this manuscript, we will analyze the case of vector fields with singularities and exploit further geometrical consequences of (\ref{integral-formula}).

\end{document}